\documentclass[a4paper,12pt]{amsart}

\title{}
\author{}

\usepackage{amsmath,amssymb,amsthm,amscd,amsfonts}
\usepackage{latexsym}

\usepackage{amsxtra}
\usepackage{eucal}

\newtheorem{defin}{Definition}
\newtheorem{prop}{Proposition}

\newtheorem*{rmk}{Remark}

\newcommand{\Nil}{\mathcal N}

\DeclareMathOperator*{\fprod}{\times}
\author{Dmitry Korb}
\title[Tadpole quiver]
{Irreducible components of the nilpotent cone of the tadpole quiver}
\begin{document}

\begin{abstract}
We exhibit the irreducible components of the nilpotent cone of the tadpole quiver of arbitrarily high codimension and introduce an approach that could allow to describe all of them.
\end{abstract}
\maketitle

\section*{Introduction}

This note deals with the nilpotent cone of the tadpole quiver, which is one of the simplest quivers with loops. The goal is to find a description for the irreducible components of the nilpotent cone. According to G.~Lusztig, the
nilpotent cone of a quiver without loops is equidimensional (lagrangian).
According to A.~Khoroshkin (unpublished), the nilpotent cone of the tadpole
quiver has some irreducible components of dimension less than expected. 
In this note we construct a series of examples of irreducible components of
arbitrarily high codimension.

The definitions of the nilpotent cone and tadpole quiver are given in Section \ref{def}. Section \ref{nil} deals with a useful property of the equations defining the nilpotent cone, which allows to describe it as a certain fiber product of relatively simple manifolds. This description, along with a stratification arising from Jordan stratification on the space of nilpotent matrices, is given in the next section and is shown to be useful in the study of irreducible components. Section \ref{num} describes an upper bound on the number of components, effectively reducing the problem to the description of a smaller variety. While the complete description is unknown, this result allows one to construct some of the irreducible components of this nilpotent cone. Two different cases corresponding to the highest and the lowest elements of Jordan stratification are considered in Section~\ref{ex}, giving the components of maximal and minimal possible codimension respectively. In particular, an example of the components of arbitrarily high codimension is described.

\section{Definitions}
\label{def}
We recall some material of~\cite{Lus}.
Let $\Gamma$ be a finite graph, with $I$ being the set of vertices and $E$ being the set of edges. Let $H$ be the set consisting of edges along with a choice of orientation on them, then an orientation of $\Gamma$ is a subset $\Omega \subset H$, such that $\overline{\Omega} \cap \Omega = \emptyset$ and $\overline{\Omega} \cup \Omega = H$, where the bar denotes the inverted orientation. The choice of orientation is equivalent to the choice of function $\omega:H \to \{\pm1\}$, such that $\omega(\overline h)=-\omega(h)$, the correspondence given as follows: $\Omega = \{h \in H|\omega(h)=1\}$. For $h \in H$, let $in(h)$ and $out(h)$ denote its incoming and outcoming vertices respectively.

The quiver is a finite graph $\Gamma$ with an orientation $\Omega$. Most of the constructions below do not, in fact, depend on the choice of orientation, so we will often refer to a graph as "quiver" without specifying the orientation, thus meaning that any orientation will do. Particular choice of orientation will be required mostly for notational purposes.

\begin{defin}
The tadpole quiver $T_n$ is a quiver with $n$ vertices $1$, $2$, $\dots$, $n$ and $n$ edges $h_1$, $h_2$, $\dots$, $h_n$, with orientation chosen, such that $in(h_i)=i$, $out(h_i)=i+1$ $\forall 1 \leq i \leq n-1$ and $in(h_n)=out(h_n)=n$.
\end{defin}

In other words, it's a quiver of type $A_n$ with an additional loop at the last vertex.

For arbitrary quiver ($\Gamma$, $\Omega$), fix a dimension vector $v \in \mathbb Z_{\geq 0}^I$ and a set of vector spaces $\{V_i|i \in I, dim(V_i)=v_i\}$. Consider the space $M_{\Gamma, v}=\bigoplus_{h \in H} Hom(V_{in(h)}, V_{out(h)})$. For an arbitrary element $B \in M_{\Gamma, v}$, $B=\bigoplus_{h \in H} B_{h}$.

\begin{defin}
An element $B \in M_{\Gamma, v}$ is called nilpotent, if $\exists N \in \mathbb N$, such that for any path $h_1 h_2 \dots h_k$, $out(h_i)=in(h_{i+1})$, $k>N$, the composition $B_{h_k} B_{h_{k-1}} \dots B_1=0$.
\end{defin}

\begin{defin}
The moment map $\mu: M_{\Gamma, v} \to \bigoplus_{i \in I} \mathfrak {gl}(V_i)$ is given by the following formula:

$$\mu(B)_i = \sum_{h \in H: in(h)=i} \omega(h)B_{\overline h}B_h$$
\end{defin}

It is obvious from this formula that $\mu$ depends on the orientation. Change of orientation is equivalent, as far as the moment map is concerned, to multiplication by $-1$ in several matrix summands of $M_{\Gamma, v}$, which does not change the resulting varieties.

\begin{defin}
The subvariety

$$\mathfrak N_{\Gamma, v} = \{B \in M_{\Gamma, v}|\mu(B)=0, B\text{ is nilpotent}\} \subset M_{\Gamma, v}$$

is called nilpotent cone corresponding to the quiver $\Gamma$.
\end{defin}
\section{Nilpotency condition}
\label{nil}
It is well known that the definition of the nilpotent cone is, in fact, redundant. To see this let us first consider a quiver of type $A_n$ with the edges being oriented in one direction. Denoting the vertices $1$, $2$, $\dots$, $n$, and the edges $h_1$, $h_2$, $\dots$, $h_{n-1}$, such that $in(h_i)=i$, we can write the condition $\mu=0$ as a set of equations $B_{\overline {h_1}}B_{h_1}=0$, $B_{h_{n-1}}B_{\overline{h_{n-1}}}=0$, $B_{h_{i-1}}B_{\overline{h_{i-1}}}=B_{\overline {h_i}}B_{h_i}$ $\forall$ $2\leq i \leq n-1$. This means that, first, for any path containing a fragment $\dots h_i \overline{h_i} \dots$ the resulting product as in the definition of the nilpotent element coincides with the product for the path obtained by replacing this fragment by $\dots \overline{h_{i-1}} h_{i-1} \dots$, and second, this product automatically equals to $0$ for any path passing through the extreme vertices. It follows that the condition $\mu=0$ implies the nilpotency condition in this case, with $N=2(n-2)$ $\forall B \in \mu^{-1}(0)$.

Now let us consider a finite connected graph $\Gamma$ with cycles generated by loops only, such that all vertices save one are of degree less than or equal to 2. This means, equivalently, that this graph consists of several type $A$ quivers and several loops linked to one initial vertex, a special case of this being the tadpole quiver. Fix an orientation $\Omega$. Denote the initial vertex $i_0$, and denote $A = \{h \in H|in(h)=i_0, out(h) \neq i_0\}$, $L = \{h \in \Omega|in(h)=out(h)=i_0\}$. As before, let $V$ and $v$ be the set of vector spaces corresponding to the vertices and their dimension vector respectively. Consider the following operators $N_h$, $h \in A \cup L$, in $End(V_{i_0})$: $\forall h \in A$, $N_h=\omega (h)B_{\overline h}B_h$; $\forall h \in L$, $N_h=[B_{\overline h},B_h]$. The condition $\mu=0$ is the same as in the case of type $A$ quiver for all vertices save the initial one, and for the latter it can be written as $\sum_{h\in A \cup L} N_h=0$.

\begin{defin}
The set of operators $X_j \in End(V)$, $j \in J$ for any vector space V and any index set J is called nilpotent, if $\exists N \in \mathbb N$ such that any composition of $k$ operators $X_j$ $\forall k>N$ is equal to $0$.
\end{defin}

\begin{prop}
\label{Local nilpotency}
Under these assumptions, $B \in \mu^{-1}(0)$ is nilpotent iff $\{N_h|h \in A\} \cup \{B_h|h \in L \cup \overline L\}$ is nilpotent.
\end{prop}

{\itshape Proof.} To prove that $B \in \mu^{-1}(0)$ is nilpotent it suffices to prove that the composition is 0 for any sufficiently long path beginning and ending in the initial vertex. Any such path can be divided into elementary parts beginning and ending in the initial vertex and not passing through it. Due to the conditions $\mu=0$ for the vertices of quiver of type $A_n$, the composition of operators from $B$ along any such elementary path would equal to either $B_h$ for $h \in L \cup \overline L$ or $N_h^k$ for some $h \in A$, $k \in \mathbb N$. Thus the condition that $\{N_h|h \in A\} \cup \{B_h|h \in L \cup \overline L\}$ is nilpotent implies $B \in \mu^{-1}(0)$ is nilpotent. The second part of the statement is straightforward. \qed

Proposition \ref{Local nilpotency} states that the nilpotency condition is essentially local, that is, the only restrictions concern adjacent vertices. In effect, it allows to construct the nilpotent cone for such quivers as a fiber product of a set of relatively simple varieties corresponding to the edges of a quiver over another set of varieties corresponding to the vertices.
\section{Decomposition}
\label{decom}
Let $U$ be a vector space. Denote $\Nil(U)$ the variety of nilpotent endomorphisms of $U$. Let $U_1$, $U_2$ be two vector spaces, denote $\mathfrak H(U_1, U_2) = \{(h, \overline h)|h \in Hom(U_1, U_2), \overline h \in Hom(U_2, U_1), \overline hh \in \Nil(U_1), h\overline h \in \Nil(U_2)\}$. Note that there exist natural morphisms $p_i$: $\mathfrak H(U_1, U_2) \to \Nil(U_i)$, $i=1,2$. Finally, for a vector space $U$, denote $\mathfrak X(U)=\{(h, \overline h)|h \in End(U), \overline h \in End(U)$, the pair $(h, \overline h)$ is nilpotent$\}$. There is a morphism $\phi: \mathfrak X(U) \to \Nil(U)$, $\phi(h, \overline h)=[h, \overline h]$.

\begin{prop}
\label{Structure of nc}
The variety $\mathfrak N_{T_n, v}$ is isomorphic to the fiber product $\mathfrak H(0, V_1) \fprod \limits_{\Nil(V_1)} \mathfrak H(V_1, V_2) \fprod \limits_{\Nil(V_2)} \dots \fprod \limits_{\Nil(V_{n-1})} \mathfrak H(V_{n-1}, V_n) \fprod \limits_{\Nil(V_n)} \mathfrak X(V_n)$.
\end{prop}

{\itshape Proof.} This is a direct corollary of Proposition \ref{Local nilpotency} and the condition $\mu=0$. For any $B \in \mathfrak N_{T_n, v}$ let $A_i=B_{\overline {h_i}}B_{h_i} \in End(V_i)$ $\forall i<n$, $A_n=[B_{h_n}, B_{\overline {h_n}}]$. Then $A_i \in \Nil(V_i)$, $A_1=0$, $B_{h_i}B_{\overline {h_i}}=A_{i+1}$ $\forall i<n$, and the triple $(B_{h_n}, B_{\overline {h_n}}, A_n)$ is nilpotent, hence $\mathfrak N_{T_n, v} \subset \mathfrak H(0, V_1) \fprod \limits_{\Nil(V_1)} \mathfrak H(V_1, V_2) \fprod \limits_{\Nil(V_2)} \dots \fprod \limits_{\Nil(V_{n-1})} \mathfrak H(V_{n-1}, V_n) \fprod \limits_{\Nil(V_n)} \mathfrak X(V_n)$. The opposite inclusion is again straightforward. \qed

\begin{rmk}
{\em
Note that similar decomposition exists for any quiver of the type considered in the previous section. One important example is $A_n$. The same arguments as before show that $\mathfrak N_{A_n, v}$ is isomorphic to the product $\mathfrak H(0, V_1) \fprod \limits_{\Nil(V_1)} \mathfrak H(V_1, V_2) \fprod \limits_{\Nil(V_2)} \dots \fprod \limits_{\Nil(V_{n-1})} \mathfrak H(V_{n-1}, V_n) \fprod \limits_{\Nil(V_n)} \mathfrak H(V_n, 0)$.
}
\end{rmk}

For a vector space $U$ of dimension $d$ consider the variety $\Nil(U)$. It admits stratification by subvarieties $\Nil_{\lambda}(U) \subset \Nil(U)$, $\forall \lambda \vdash d$, where $\Nil_{\lambda}(U)$ is a locally closed subset of all nilpotent operators $A$ with Jordan form corresponding to the partition $\lambda$, that is, the sizes of Jordan blocks of $A$ are equal to the elements of partition or, equivalently, $\dim(Ker(A^i))-\dim(Ker(A^{i-1}))=(\lambda^T)_i$, where $\lambda^T$ is a transposed partition.

The map $\phi: \mathfrak X(U) \to \Nil(U)$ described earlier allows to introduce similar stratification $\mathfrak X_{\lambda}(U) = \phi^{-1}(\Nil_{\lambda}(U)) \subset \mathfrak X(U)$. Note that $\phi$ is not surjective. In particular, $\mathfrak X_{\lambda}(U) = \emptyset$ iff $\lambda$ has more than $\frac {d+1}2$ columns.

Now let $U_1$, $U_2$ be two vector spaces of finite dimensions $d_1$, $d_2$ respectively. The stratifications of the varieties $\Nil(U_i)$ along with two natural morphisms $p_i$ lead to a stratification $\mathfrak H_{\lambda, \mu}(U_1, U_2) = p_1^{-1}(\Nil_{\lambda}(U_1)) \cap p_2^{-1}(\Nil_{\mu}(U_2)) \subset \mathfrak H(\Nil(U_1), \Nil(U_2))$ $\forall \lambda \vdash d_1$, $\mu \vdash d_2$.

The same is applicable to the quiver varieties. In case of quiver $\Gamma$ of type $A_n$ or $T_n$, the operators $B_{\overline {h_i}}B_{h_i} \in End(V_i)$ give a map $\mathfrak N_{\Gamma, v} \to \Nil(V_i)$ and the preimages of the strata $\Nil_{\lambda^i}(V_i)$ form a stratification $\mathfrak N_{\Gamma, \Lambda} \subset \mathfrak N_{\Gamma, v}$ itself, for $\Lambda=(\lambda^1 \vdash v_1, \lambda^2 \vdash v_2, \dots, \lambda^{n-1} \vdash v_{n-1}, \lambda^n \vdash v_n)$. For most $\Lambda$, however, the corresponding stratum turns out to be empty. For example, as $B_{\overline {h_1}}B_{h_1}=0$, the stratum is empty if $\lambda^1 \neq (1^n)$, and for $A_n$ the same holds for $\lambda^n$. Note that the description of the nilpotent cones as fiber products is compatible with Jordan stratification, in particular, the stratum $\mathfrak N_{\Gamma, \Lambda}$ is isomorphic to the fiber product of the strata of the factors over the corresponding strata of $\Nil(V_i)$, and it is empty if and only if at least one factor is empty. Let us call such a set of partitions $\Lambda$ a \textit{multipartition} of $v$, and write $\Lambda \vdash v$.

\begin{defin}
If $\lambda \vdash d_1$ and $\mu \vdash d_2$ are two partitions, a proper pairing between $\lambda$ and $\mu$ is a correspondence between some of the elements of the partitions, such that any element $\lambda_i \geq 2$ of a partition, corresponds to an element $\mu_j$ of another partition, such that $|\lambda_i-\mu_j| \leq 1$, and elements equal to $1$ may not have a pair.
\end{defin}

The geometric significance of this notion is as follows: $\mathfrak H_{\lambda, \mu}(U_1, U_2) \neq \emptyset$ if and only if there exists a proper pairing between $\lambda$ and $\mu$. This follows easily from the classification of irreducible nilpotent representations of Jordan quiver. This implies that $\mathfrak N_{A_n, \Lambda} \neq \emptyset$ if and only if $\forall i$ there exists a proper pairing between $\lambda^i$ and $\lambda^{i+1}$. Let us call a multipartition satisfying this condition \textit{admissible}. 
In particular, for any admissible $\Lambda$, $\lambda^1$ and $\lambda^n$ should contain only one column, as it is the only partition admitting a proper pairing with an empty partition of $0$.

Let us denote $I_{\lambda, \mu}^{\mathfrak H}(U_1, U_2)$ the set of irreducible components of $\mathfrak H_{\lambda, \mu}(U_1, U_2)$, denote $I_{\lambda}^{\mathfrak X}(U)$ the set of irreducible components of $\mathfrak X_{\lambda}(U)$ and denote $I_{\Gamma, \Lambda}^{\mathfrak N}$ the set of irreducible components of $\mathfrak N_{\Gamma, \Lambda}$. Note that $I_{A_n, \Lambda}^{\mathfrak N}$ is in natural bijection with $\prod_{i=1}^{n-1} I^{\mathfrak H}_{\lambda^i, \lambda^{i+1}}(V_i, V_{i+1})$, and the dimension of $\gamma \in I_{A_n, \Lambda}^{\mathfrak N}$ is equal to $\sum_{i=1}^{n-1} \dim  \gamma_i - \sum_{i=1}^n \dim \Nil_{\lambda^i}(V_i)$, where $\gamma_i \in I^{\mathfrak H}_{\lambda^i, \lambda^{i+1}}(V_i, V_{i+1})$ are the corresponding irreducible components of the factors. In the same way,  $I_{T_n, \Lambda}^{\mathfrak N}=\left (\prod_{i=1}^{n-1} I^{\mathfrak H}_{\lambda^i, \lambda^{i+1}}(V_i, V_{i+1}) \right ) \times I_{\lambda^n}^{\mathfrak X}(V_n)$, and the dimension of every component is equal to the sum of dimensions of the corresponding components of factors minus $\sum_{i=1}^n \dim \Nil_{\lambda^i}(V_i)$.

Note that every irreducible component of $\mathfrak N_{\Gamma, v}$ is a closure of an irreducible component of one of $\mathfrak N_{\Gamma, \Lambda}$, and every irreducible component of $\mathfrak N_{\Gamma, \Lambda}$ belongs to an irreducible component of $\mathfrak N_{\Gamma, v}$. It allows one to derive the number and dimensions of irreducible components of $\mathfrak H_{\lambda, \mu}(U_1, U_2)$ from $A_n$ type quivers, and could lead to the description of the number and dimensions of irreducible components of $\mathfrak N_{T_n, v}$ if one can find them for $\mathfrak X_{\lambda}(U)$.

\begin{prop}
\label{Elementary}
$\forall \mu \vdash d_1$, $\eta \vdash d_2$ the irreducible components of $\mathfrak H_{\mu, \eta}(U_1, U_2)$ have the same dimension $d_1d_2+\frac12 \dim\Nil_{\mu}(U_1)+\frac12 \dim\Nil_{\eta}(U_2)$. The number of these components is denoted $\chi(\mu, \eta)$.
\end{prop}

{\itshape Proof.}
If $\lambda$ is a partition, let $t(\lambda)$ be a partition obtained by omitting the first column of $\lambda$. Note that there exists a proper pairing between $\mu$ and $\eta$ if and only if $t(\mu) \subset \eta$ and $t(\eta) \subset \mu$ as Young diagrams.

Recall that all irreducible components of $\mathfrak N_{A_n, v}$ have the same dimension equal to $\sum_{i=1}^{n-1}v_iv_{i+1}$.

If $\mu$ and $\eta$ are two partitions with a proper pairing, consider a quiver $A_n$ for $n=\mu_1+\eta_1$, dimension vector $v$ and a multipartition $\Lambda \vdash v$, such that $\lambda^i=t^{\mu_1-i}(\lambda)$ for $i \leq \mu_1$ and $\lambda^i=t^{i-\mu_1-1}(\eta)$ for $i > \mu_1$. For such data $\mathfrak N_{A_n, \Lambda}$ is a maximal nonempty element of Jordan filtration on $\mathfrak N_{A_n, v}$, in particular, it is nonempty open subset of the latter, and the closures of its irreducible components are irreducible components of $\mathfrak N_{A_n, v}$, hence they are equidimensional. Thus, as $I_{A_n, \Lambda}^{\mathfrak N}=\prod_{i=1}^{n-1} I^{\mathfrak H}_{\lambda^i, \lambda^{i+1}}(V_i, V_{i+1})$, irreducible components of $\mathfrak H_{\mu, \eta}(V_{\mu_1}, V_{\mu_1+1})$ are also equidimensional. The dimension is easily obtained by induction, starting from $\dim \mathfrak H_{(1^{d_1}),(1^{d_2})}(V_1, V_2)=d_1d_2$. \qed

\section{Number of components}
\label{num}
Let us modify the notion of admissible multipartition to describe irreducible components of $\mathfrak N_{T_n, v}$. A multipartition is called $(\mu, \eta)$\textit{-admissible}, if $\lambda^1=\mu, \lambda^n=\eta$ and $\forall i$, $1 \leq i \leq n-1$, there exists a proper pairing between $\lambda^i$ and $\lambda^{i+1}$. Let $\mathcal A_v(\mu, \eta)$ be a set of all $(\mu, \eta)$-admissible $\Lambda \vdash v$. Let $\chi(\Lambda)=\prod_{i=1}^{n-1}\chi(\lambda^i, \lambda^{i+1})$, for $\chi(\lambda, \mu)$ defined in Proposition \ref{Elementary}. Finally, let $\psi(\lambda)=\left| I_{\lambda}^{\mathfrak X}(U)\right|$ for a vector space $U$ of dimension $|\lambda|$.

The stratification introduced in the previous section along with representation of the nilpotent cone as fiber product gives a way to compute the number of irreducible components of the nilpotent cone $\mathfrak N_{\Gamma, v}$ as well as their dimensions.

If $\Gamma = A_n$, there is $K(v)=\sum_{\Lambda \in \mathcal A_v((1^{v_1}), (1^{v_n}))} \chi(\Lambda)$ components, their dimensions are all equal to $\sum_{i=1}^{n-1}v_iv_{i+1}$. On the other hand, it is well known that $K(v)$ is a Kostant partition function. It gives one a way to compute the values of $\chi$: for any $\lambda$ and $\mu$, consider a quiver as in the proof of Proposition \ref{Elementary}. For pairs of partitions $(\lambda, \mu)$ and $(\eta, \nu)$, let $(\lambda, \mu) \geq (\eta, \nu)$ if $\lambda \succeq \eta$ and $\mu \succeq \nu$, where $\lambda \succeq \eta$ denotes usual dominance ordering if $|\lambda|=|\eta|$, and denotes $|\lambda| > |\eta|$ otherwise. $K(v)=\sum_{\Lambda \in \mathcal A_v((1^{v_1}), (1^{v_n}))} \chi(\Lambda)$, and for all $\Lambda$ encountered in this equation $\forall i$, $(\lambda^i, \lambda^{i+1}) \leq (\lambda, \mu)$. Thus, if we know $\chi(\eta, \nu)$ for every pair $(\eta, \nu) < (\lambda, \mu)$, we can derive $\chi(\lambda, \mu)$ from this equation. There are finitely many such pairs $(\eta, \nu)$, so this allows to compute these numbers recursively. For example, $\chi((1^a), (1^b))=K((a, b))= \min (a, b)+1$; considering the quiver $A_3$ one can find that $\chi((2^a1^b), (2^c1^d))=0$ if $a+b<c$ or $c+d<a$, that is, when there is no proper pairing, and equals to $\min(a, b)+\min(a+b-c, c+d-a)+1$ otherwise.

If $\Gamma = T_n$, the number of components of a stratum $I_{\Gamma, \Lambda}^{\mathfrak N}$ is equal to $\chi(\Lambda)\psi(\lambda^n)$, and the dimension of any such component is equal to the dimension of underlying component of $\mathfrak X_{\lambda^n}(V_n)$ plus $\left(\sum_{i=1}^{n-1}v_iv_{i+1}\right)-\frac12\dim\Nil_{\lambda^n}(V_n)$. Any irreducible component of $\mathfrak N_{T_n, v}$ coincides with a closure of an irreducible component of some stratum, so this effectively gives an upper bound on the set of irreducible components, because an irreducible component of a stratum may belong to the closure of irreducible component of a higher stratum. In low-dimension cases it does not in fact occur: the irreducible components of the higher strata can be shown to have dimension not greater than that of the components of a lower stratum, that is, $\frac12\dim\Nil_{\lambda^n}(V_n)$ increases faster than the dimensions of the components of $\mathfrak X_{\lambda^n}(V_n)$. It is not known if it holds true in general, and the set $I_{\lambda^n}^{\mathfrak X}(V_n)$ is also unknown.

\section{Examples}
\label{ex}
Note that the variety of pairs of commuting nilpotent matrices $\mathfrak X_{(1^d)}(U)$ is irreducible of dimension $d^2-1$. There is a proof of this fact in \cite{Bar}.

Let $v_n$ be equal to 1 or 2. In this case $\phi \equiv 0$, the variety $\mathfrak X(V_n)=\mathfrak X_{(1^{v_n})}(V_n)$ is irreducible of dimension $v_n^2-1$, and $\mathfrak N_{T_n, v} \simeq \mathfrak N_{A_n, v} \times \mathfrak X(V_n)$. In particular, all irreducible components have the same dimension. If $v_n \geq 3$, $\mathfrak X(V_n)$ is no longer equal to $\mathfrak X_{(1^{v_n})}(V_n)$, but there exists closed embedding $\mathfrak N_{A_n, v} \times \mathfrak X_{(1^{v_n})}(V_n) \hookrightarrow \mathfrak N_{T_n, v}$. The dimension of image of this embedding coincides with the dimension of the nilpotent cone, so one obtains a set of equidimensional irreducible components corresponding to irreducible components of $\mathfrak N_{A_n, v}$.

Let $\mathfrak A_n(v) = \mathfrak H(0, V_1) \fprod \limits_{\Nil(V_1)} \mathfrak H(V_1, V_2) \fprod \limits_{\Nil(V_2)} \dots \fprod \limits_{\Nil(V_{n-1})} \mathfrak H(V_{n-1}, V_n)$, and $\mathfrak A_{n, \lambda}(v) \subset \mathfrak A_n(v)$ is a stratification obtained from Jordan stratification of $\Nil(V_n)$.

Consider the case of tadpole quiver with n vertices. Let $v_i \geq 2i$, $v_n=2n$ or $v_i \geq 2i-1$, $v_n=2n-1$, $\lambda=(\lfloor \frac{v_n+1}2 \rfloor, \lfloor \frac{v_n}2 \rfloor)$. Note that $\mathfrak X_{\lambda}(v_n)$ is an open dense subset of $\mathfrak X(v_n)$, in particular, $\dim \mathfrak X_{\lambda}(v_n)=\frac32v_n(v_n-1)$. Note also that $\dim \Nil_{\lambda}(v_n)=4\lfloor \frac{v_n}2 \rfloor \lfloor \frac{v_n-1}2 \rfloor$. Hence $\mathfrak A_{n, \lambda}(v) \fprod \limits_{\Nil_{\lambda}(v_n)} \mathfrak X_{\lambda}(v_n)$ is an open subset in $\mathfrak N_{T_n, v}$ of dimension $\sum_{i=1}^{n-1} v_i v_{i+1}+v_n^2-n$. It is nonempty due to the conditions on $v$, and its irreducible components are irreducible components of $\mathfrak N_{T_n, v}$ of codimension $n-1$. The number of these components is equal to $\sum_{\Lambda \in \mathcal A_v((1^{v_1}), \lambda)} \chi(\Lambda)$.

\bigskip

\footnotesize{ National Research
University Higher School of Economics,\\ Department of
Mathematics, 7 Vavilova st, Moscow 117312 Russia;\\ 
{\tt korbdv@gmail.com}}

\end{document}